\documentclass[a4paper,12pt,notitlepage]{article}

\usepackage{amsmath,amssymb,amsthm,amscd}
\usepackage{oldgerm} 
\usepackage{epsfig}
\usepackage[all]{xy}
\usepackage{a4wide}
\usepackage{xspace} 
\usepackage{makeidx}
\usepackage{shapepar}

\CompileMatrices 

\newtheorem{thm}{Theorem}[section]
\newtheorem{cor}[thm]{Corollary}
\newtheorem{lem}[thm]{Lemma}

\theoremstyle{definition}
\newtheorem{df}[thm]{Definition}
\theoremstyle{remark}
\newtheorem{rem}[thm]{Remark}

\newtheorem{ex}[thm]{Example}

\newcommand{\Aff}{\operatorname{Aff}}

\newcommand{\Id}{\operatorname{Id}}

\newcommand{\bZ}{{\mathbb Z}}
\newcommand{\bH}{{\mathbb H}}
\newcommand{\bR}{{\mathbb R}}

\newcommand{\bN}{{\mathbb N}}

\newcommand{\inv}{{\operatorname{Inv}}}

\newcommand{\Cucb}{\operatorname{C_{ucb}}}

\newcommand{\glnr}{\mbox{$GL(n, {\mathbb R})$}}
\newcommand{\lp}{\mbox{$\bigoplus^\infty_{n=1} L^{2n}(G, \mu)$ (in the $l^2$ sense)}}

\newcommand{\plig}{{\bf plig} metric\xspace}
\newcommand{\pligs}{{\bf plig} metrics\xspace}
\newcommand{\lcsc}{ locally compact, second countable\xspace}

\newcommand{\expgrow}{the $d$-balls have exponentially controlled growth\xspace}

\newcommand{\bref}[1]{(\ref{#1})}

\title{Proper metrics on locally compact groups, and proper affine isometric
  actions on Banach spaces.}
\author{  Uffe Haagerup and Agata Przybyszewska\thanks{Supported by the
  Danish National Research Foundation}}
\date{}

\begin{document}

\maketitle



\begin{center}
{\bf Abstract}

\vspace{0.7cm}

\parbox{0.8\linewidth}{\small
 In this article it is proved, that every
  \lcsc group has a left invariant metric $d$, which generates the topology
  on $G$, and which is proper, ie. every closed $d$-bounded set in $G$ is
  compact. Moreover, we obtain the following extension of a result due
  to N. Brown and E. Guentner \cite{brown_guentner}:
     Every \lcsc  $G$ admits a proper affine action on the reflexive and
    strictly convex Banach space  
 $$\bigoplus^{\infty}_{n=1} L^{2n}(G, d\mu),$$
where the direct sum is taken in the $l^2$-sense.
}
\end{center}

\section{Introduction}
We consider a special class of metrics on  second countable,
locally compact groups, namely proper left invariant metrics which generate
the given topology on $G$, and we will denote such a metric a \plig. 

It is fairly easy to show, that if a locally compact group $G$ admits a
\plig, then $G$ is second countable. Moreover, any two \pligs $d_1$ and
$d_2$ on $G$ are coarsely equivalent, ie. the identity map on $G$ defines a
coarse equivalence between $(G, d_1)$ and $(G, d_2)$ in the sense of
\cite{roe}.

In \cite[pp. 14-16]{lubotzky_etal} Lubotzky, Moser and Ragunathan shows, that every compactly
generated second countable group has a \plig. 
Moreover, Tu shoved in \cite[lemma 2.1]{tu}, that every contable discrete
group has a \plig.

The main result of this paper
is that every \lcsc group $G$ admits a \plig. Moreover a \plig $d$ can be
chosen, such that the $d$-balls have exponentially controlled growth in the
sense that 
$$ \mu( B_d(e,n)) \leq \beta\cdot e^{\alpha n}, \quad n \in \bN,
$$
for suitable constants $\alpha$ and $\beta$. Here $\mu$ denotes the Haar measure on $G$.

In \cite{brown_guentner}, Brown and Guentner proved that for every contable
discrete group $\Gamma$ there exists a sequence of numbers $p_n\in (1,
\infty)$ converging to $\infty$ for $n\rightarrow\infty$, such that
$\Gamma$ has a proper affine action on the reflexive and strictly convex
Banach space 
$$ X_0 = \bigoplus^\infty_{n=1} l^{p_n}(\Gamma), $$
where the direct sum is in the $l^2$-sense. 

By similar methods, we prove that every second countable group has a proper
affine action on the reflexive and strictly convex Banach space
$$ X = \bigoplus^\infty_{n=1}\quad L^{2n}(G, d\mu),$$
where the sum is taken in the $l^2$-sense. However, in order to prove, that the exponents $(p_n)^\infty_{n=1}$ can be
chosen to be $p_n=2n$, it is essential to work with a \plig on $G$ for
which the $d$-balls have exponentially controlled growth. As a corollary we obtain, that a second countable locally compact group $G$ has a uniform embedding in the above Banach space $X$

Note, that the Banach spaces $X$ and $X_0$ above are not uniformly
convex. Kasparov and Yu have recently proved, that the Novikov conjecture holds for every discrete countable group, which has a uniform embedding in a uniformly convex Banach space (cf. \cite{kasparov_yu}).

\subsection*{Acknowledgments}
We would like to thank Claire Anantharaman-Delaroche, George
Skandalis, Alain Valette and  Guoliang Yu  for stimulating
mathematical conversations.

\section{Coarse geometry and \pligs}

First, let us fix notation and basic definitions:

\begin{df} Let $G$ be a topological group, ie. a group equipped with a
  Hausdorff topology, such that the map $ (x,y) \rightarrow x\cdot y^{-1}$  is
  continuous. If $G$ is equipped with a metric $d$, we put
$ B_d(x, R) = \{ y\in G: d(x,y)<R\}$ and  $D_d(x, R) = \{ y\in G: d(x,y)\leq R\}$.

  \begin{itemize}
  \item $G$ is {\em locally compact}\index{locally compact} if  every $x\in G$ has a
    relative compact neighbourhood.
    

    \item We say that the metric $d$ {\em induces the topology}\index{metric!induces the topology} of $G$ if
    the topology generated by the metric $\tau_d$ coincides with the original
    topology $\tau$.

    
  \item  The metric $d$ is said to be {\em left invariant}\index{metric!left
      invariant} if for all $g, x, y \in G$
    we have that
    $$
    d(x, y) = d(g\cdot x, g\cdot y).$$

   \item  Following \cite[p. 10]{roe}, a metric space is called {\em
    proper}\index{space!proper} if all closed bounded sets are compact. When $G$ is a group, 
     this reduces by the left invariance of the group metric  to the requirement, 
     that for every $M>0$ all the closed balls  
$$   D(e, M) = \{ h\in G: d(e, h)\leq M \}$$
 are compact.
 \end{itemize}
   
   \end{df}

We will work with a special class of metrics on \lcsc groups, which we
define here:
\begin{df}
  Let $G$ be a topological group. A \plig $d$ on $G$ is a metric on $G$,
  which is  {\bf
        P}roper, {\bf L}eft {\bf I}nvariant and {\bf G}enerates the
      topology\index{metric!plig}\index{plig \see{metric!plig}}.
\end{df}

M. Gromov started investigating asymptotic invariants of groups,
particularly fundamental groups of manifolds.

This lead to the development of coarse geometry - or large scale geometry
\index{geometry!large scale \see{coarse!geometry}}. Coarse geometry  studies global properties of metric spaces, neglecting
small (bounded) variations of these spaces. The properties and invariants
in coarse geometry are treated in the limit at infinity, as opposed to the
traditional world of topology, which focuses on the local structure of the
space.

\begin{df}[\cite{roe}]
  Let $(X, d_X)$ and $(Y, d_Y)$ be metric spaces.
\begin{itemize}
\item A map $f:X\rightarrow Y$ is called {\em
    uniformly expansive}\index{map!uniformly expansive} if
  $$ \forall_{ R>0} \exists_{S>0} \quad\mbox{ such that }\quad
  d_X(x,x')\leq R\Rightarrow d_Y(f(x), f(x'))\leq S.$$

\item A map $f:X\rightarrow Y$ is called {\em
    metrically proper}\index{map!metrically proper} if
  $$ \forall_{B\subset Y} \qquad B\mbox{ is bounded }\Rightarrow  f^{-1}(B)\subset
X \mbox{ is bounded}.$$

\item  A map $f:X\rightarrow Y$ will by called a {\em coarse map}\index{map!coarse} if it is
  both metrically proper and uniformly expansive.

\item Two coarse maps $h_0, h_1:X\rightarrow Y$ are said to be {\em coarsely
  equivalent}\index{map!coarsely equivalent} when
$$\exists_{C>0}\forall_{x\in X}:\quad d_Y(h_0(x), h_1(x))<C. $$
We denote the relation of coarse equivalence by '$\sim_c$'.

\item
The spaces $X$ and $Y$ are said to be {\em coarsely
  equivalent}\index{space!coarsely equivalent} if there
exist coarse maps  $f:X\rightarrow Y$, and $ g:Y\rightarrow X$ such that
$$ f\circ g \sim_c
\Id_Y\mbox{ and } g\circ f\sim_c\Id_X.$$

\item    The {\em coarse structure}\index{coarse!structure} of  $X$ means the coarse
equivalence class 
of the given metric, $[d_X]_{\sim}$. 

  \end{itemize}
\end{df}

\begin{ex}
  It is well known  that 
  \begin{equation*}
    (0,1) \sim_h \bR\quad\mbox{ and }\quad
  \bZ \not\sim_h\bR,
\end{equation*}
where $\sim_h$ means that the two spaces in question are homeomorphic.

It is easy to see, that in the coarse case we have:
 \begin{equation*}
    (0,1) \not\sim_c \bR\quad\mbox{ and }\quad
  \bZ \sim_c\bR.
\end{equation*}
\end{ex}

 The reason for working with coarse geometry is, that many geometric group
 theory properties are coarse invariants.  Examples of coarse invariants are: property A \cite{higson_and_roe},
  asymptotic dimension \cite{bell_drani_asdim} and 
  change of generators for a finitely generated group.

M. Gromov has suggested in \cite{gromov_novikov} to solve the Novikov
conjecture by considering classes of groups admitting uniform embeddings
into Banach spaces with various restraints. 

G. Yu proved in \cite{yu:uniform_embedding}, that the Novikov Conjecture is
true for a space that admits a uniform embedding into a Hilbert space. 
This result was strenghtened in \cite{kasparov_yu}, where it is shown that
the Novikov Conjecture is
true for a space that admits a uniform embedding into a uniformly convex
Banach space.

Together with the fact, that exact groups admit a uniform embedding into a
Hilbert space, see \cite{kaminker_and_guentner}, makes  uniform embedding  extremely
interesting to study.

The idea of a uniform embedding is to map a metric space  $(X, d_X)$ into
a metric space
$(Y, d_Y)$ in such a
way that the large-scale geometry of $X$ is preserved.

This means for instance that we are not
allowed to ``squease'' unbounded segments into a point, and we are not
allowed to ``blow up'' bounded segments to unbounded - the limits at infinity must be preserved.

\begin{df}
A map $f:X\rightarrow Y$ will by called a {\em uniform
    embedding}\index{uniform embedding} if
  there exist non-decreasing
     functions $\rho_-,
  \rho_+:[0,\infty)\rightarrow [0, \infty)$ such that
$$ \lim_{t\rightarrow\infty} \rho_-(t) = \lim_{t\rightarrow\infty}\rho_+(t) = \infty $$
and
  \begin{equation}
    \label{eq:df_ue_1}
    \rho_-(d_X(x,y)) \leq d_Y(f(x), f(y)) \leq \rho_+(d_X(x,y))
  \end{equation}

   When   $f:X\rightarrow Y$ is a coarse map, we will denote a map
  $\phi:f(X)\rightarrow X$ a section of $f$
  if it fullfills that  
$$ f\circ \phi = id_{f(x)}. $$
  The set of sections of the map $f$ will be denoted by $\inv(f)$. 
\end{df}

\begin{ex}
   It is easy to see that we can not
have a uniform embedding of the free group $F_2$ in any $\bR^n$. On the
other hand it was
shown in \cite{haagerup}, \cite[p. 63]{harpe_valette} that $F_2$ has a uniform embedding
in the infinite dimensional Hilbert space $\bH$. 


\end{ex}

The following is a folklore
lemma, as different definitions of uniform embedding are used in the
litterature, see \cite{gromov}, \cite{higson_and_roe} and
\cite{kaminker_and_guentner}, see \cite{przybysz_thesis} for a detailed proof. 

\begin{lem}\label{l:coarse_eq}
   Let $(X, d_X)$ and $(Y, d_Y)$ be metric spaces, and $f:X\rightarrow Y$ a map.
   The following are equivalent:
   \begin{enumerate}

   \item $f:X\rightarrow Y$ is a uniform embedding.

   \item (Guentner and Kaminker version) 
     $f$  is uniformly expansive, and 
\begin{equation}\label{eq:kaminker}
     \forall_{S>0}\exists_{R>0}\quad d_X(x,y)\geq R\Rightarrow d_Y(f(x),
       f(y))\geq S.
     \end{equation}
     
   \item (Higson and Roe version)
    The map $f$ is uniformly expansive, and so is any $\phi\in\inv(f)$.  

   \item $f$ is a coarse equivalence between  $X$ and $f(X)$ and any
     section $\phi\in\inv(f)$ is a coarse map.
   \end{enumerate}
\end{lem}

\begin{thm}\label{thm:plig_coarse_eq}
   Let $G$ be a locally compact, second countable group.
   Assume that $d_1$ and $d_2$ are \pligs on $G$. 

   Then the identity map
   $$Id:(G, d_1)\rightarrow (G,d_2)$$
   is a coarse equivalence.
 \end{thm}

 \begin{proof}


   To establish the coarse equivalence of the metric spaces, it is by
   the 3rd case of lemma~\bref{l:coarse_eq}, enough to show that $ \Id:(G,
   d_1)\rightarrow (G, d_2)$ and $ \Id:(G, d_2)\rightarrow (G, d_1)$ are
   both uniformly expansive, that is:
   \begin{eqnarray*}
      \forall_{ R>0} \exists_{S>0} \qquad
  d_1(x,y)\leq R\Rightarrow d_2(x, y)\leq S.\\
    \forall_{ R>0} \exists_{S>0} \qquad
  d_2(x,y)\leq R\Rightarrow d_1(x, y)\leq S.\\
   \end{eqnarray*}

Since both $d_1$ and $d_2$ generates the topology of $G$, the identity map 
  $$Id:(G, d_1)\rightarrow (G,d_2)$$
is a homeomorphism, and the maps $\phi_i: (G, d_i) \rightarrow \bR_+$ given
by 
$$ \phi_i(x) = d_i(e, x), \quad x\in G $$
are continuous.

Let $R>0$. Since the closed $d_1$-ball 
$$ D_1(x, R) =  \{ x\in G: d_1(x,e)\leq R  \}    $$
is compact, $\phi_2$ attains a maximum values $S$ on $D_1(x, R)$. Moreover
\begin{equation*}
   d_1( e, x)\leq R \Rightarrow d_2(e,x)\leq  S.
\end{equation*}

Hence, by the left invariance of $d_1$ and $d_2$ 
\begin{equation*}
  d_1(x,y) \leq R \Rightarrow d_2(x,y) \leq S, \quad \forall_{x,y\in G}.
\end{equation*}

 Since $R$ was arbitary, this shows the uniform expansiveness of $\Id$. 

    By reversing the roles of $d_1$ and $d_2$ in the last argument, we also
    obtain that the inverse map is uniformly expansive.
    
 \end{proof}

 \begin{rem}
   In the special case of a countable discrete group, Theorem
   \bref{thm:plig_coarse_eq} was proved by J. Tu (cf. the uniqueness part
   of \cite[lemma 2.1]{tu}).
 \end{rem}

\section{Bounded geometry on locally compact groups.}

The purpose of this section is to show, that a \plig  implies
bounded geometry on a \lcsc group. Let us begin by defining bounded
geometry, which is a concept from the world of coarse geometry.
\begin{df}\index{bounded geometry}
  Following \cite[p. 13]{roe}, the metric space $(G, d)$ is
    said to have {\em bounded geometry} if it is coarsely equivalent to a
    discrete space $(Q, d_Q)$, where for every
    $M>0$ there exists constants $\Gamma_M$ such that
    $$
    \forall_{q\in Q}\quad |D(q, M)| = |\{ p\in Q: d(q, p)\leq M
    \}| \leq \Gamma_M .$$
\end{df}

Note that $(G,d)$ is a finitely generated group, and $d$ a word length
  metric. Then $G$ has bounded geometry.

  \begin{df}\index{coarse!lattice}
    Let $(Y,d)$ be a metric space.
    We say that a discrete space $X\subset Y$ is a {\em coarse
  lattice}, if
\begin{equation}
  \exists_{\lambda\in\bR} \forall_{y\in Y}\exists_{x_y\in X}\quad d(x_y,y) \leq \lambda.
\end{equation}
  \end{df}

\begin{lem}\label{lem:proper_gives_bdg}
  Let $G$ be a locally compact group, and let $d$ be a \plig on $G$. Then
  $(G, d)$ is second countable and has bounded geometry.
\end{lem}

\begin{proof}

  First we observe that the conditions on $(G,d)$ imply, that $G$ is second
  countable. We can write $G$ as a union of compact metric spaces, namely:
$$G=\bigcup^\infty_{n=1}D_d(e,n), $$ 
where we have that each $D_d(e,n)$ is a compact metric space, since $d$ is
proper. 

Now, since  every compact metric space is seperable, see
\cite[Theorem 4.1.17, page. 297]{engelking}, we can conclude that 
every $D_d(e,n)$ has a countable dense subset. Hence $G$ is seperable, and
since for any metric space second countability is equivalent to seperability, see
\cite[Theorem 16.11, page 112]{willard}, it follows that
$G$ is second countable. 

We will now show that $G$ has bounded geometry by constructing a countable coarse lattice
$X\subset G$ such  that $X$ has bounded geometry.

  Let $X= \{x_i\}_{i\in I}$ be a maximal family of elements from $G$, such
  that $d(x_i, x_j)\geq 1, i\neq j$. Since $G$ is seperable,
  the index set $I$ is at most countable.

  By maximality of $X$, we have that
  $$ G = \bigcup_{i\in I} B_d(x_i,1).$$
  If  we had that $|I|<\infty$, then $G =  \cup_{i\in I}
  \overline{B_d(x_i,1)} $ would be a
  finite union of compact sets, and therefore compact, and hence
  bounded. Therefore
  $G$ is coarsely equivalent to $\{\bullet\}$ if $I$ is finite, and
 therefore $G$ has bounded geometry.

  Let us therefore assume, that $|I|=\infty$, and we can use $\bN$ instead
  of the index set $I$. Let us construct the coarse equivalences between
  $X$ and $G$:
 
Start by setting
 \begin{eqnarray*}
   A_1 &=& B_d(x_1, 1)\\
      A_2 &=& B_d(x_2, 1)\setminus A_1\\
      &\ldots&\\
   A_n &=& B_d(x_n, 1)\setminus \bigcup^{n-1}_{i=1} A_i
 \end{eqnarray*}

We have now that the family
  $\{A_n\}_{n=1}^\infty$ in $G$ fullfills the following:
  \begin{equation*}
    \begin{cases}
      A_n \cap A_m = \emptyset & \mbox{ if }m\neq n\\
      x_n\in A_n\subset B_d(x_n, 1)&\mbox{for all }n\in\bN\\
      \cup_{n=1}^\infty A_n = G\\
    \end{cases}
  \end{equation*}

 Now equip the set $X = \{x_i\}_{i\in I}$ with the metric inherited from
  $(G,d)$. Define:
  \begin{eqnarray*}
    \phi:G\rightarrow X,&\mbox{ by }& \phi(x) = x_n\quad\mbox{when }x\in A_n\\
 \psi:X\rightarrow G,&\mbox{ by }& \phi(x_n) = x_n\quad\mbox{for } n\in\bN
  \end{eqnarray*}

Remark that both $\psi$ and $\phi$ are coarse maps.  We have from the
construction of $\phi$ and $\psi$ that 

$$ \phi\circ\psi = \Id_X $$
and 
$$ \forall_{x\in G}\qquad d( \psi\circ\phi(x), x) \leq 1  $$
  and therefore we see, that the spaces $X$ and $G$ are coarsely equivalent.

Now, we have to show that the set $X$ indeed has
bounded geometry. Let $M>0$ be given, and let us look at the disks of radius $M$ in $X$: 
\begin{equation}
  D_X(q, M) = \{ x_n\in X:  d(q, x_n) \leq M, \quad n\in \bN\}.
\end{equation}
For every $M>0$ we need to find a constant $\Gamma_M$
such that 
$$ \sup_{q\in X}|D_X(q, M) | \leq \Gamma_M.$$

Since $d(x_n, x_m)\geq 1$ when $n\neq m$, the balls $B(x_n, \frac{1}{2})$
are disjoint. Moreover, we have that 
\begin{equation*}
  \bigcup_{x_n\in D_X(q, M) } B(x_n, \frac{1}{2}) \subset B(q, M+\frac{1}{2}).
\end{equation*}

Let $\mu$ denote the Haar measure on $G$, then we have that 
\begin{equation*}
  \sum_{x_n\in D_X(q, M)} \mu( B(x_n, \frac{1}{2})) \leq \mu( B(q, M+\frac{1}{2}))
\end{equation*}

Since the number of terms in the sum above is equal to $|D(q,M)|$, we get
by the left invariance of the Haar measure, that
$$ |D(q, M)|\cdot\mu(B(e, \frac{1}{2})) \leq \mu(B(e,M+\frac{1}{2})). $$
Hence
$$ \sup_{q\in X}  |D(q, M)|\leq \frac{ \mu(B(e,M+\frac{1}{2}))}{
  \mu(B(e,\frac{1}{2}))} < \infty.$$
Therefore we see that $(X, d)$ has bounded geometry, and we can conclude
that $(G,d)$ also has bounded geometry.
\end{proof}

   \begin{ex}\label{ex:prop_not_bdg} 
Remark, that lemma \bref{lem:proper_gives_bdg} is not true for
     general metric spaces, as
     we can find an example of a metric space $X$ that is proper, but does
     not have bounded geometry:

  Consider the triple $(D_n, d_n, x_n)$, where $D_n$ is the discrete space
  with $n$ points, equipped with the discrete metric:
  $$d_n(x, y) = \begin{cases} 1 & x\neq y\\
    0 & x= y\end{cases},$$
and where $x_n \in D_n$.

  Let $X = \coprod_{n\in \bN} D_n$, and equip $X$  with the following metric:
  $$ d(z, y) =
    d_{j(z)}(z, x_{j(z)})+ |j(z)-j(y)|+ d_{j(y)}(y, x_{j(y)}) \mbox{ where } j(x) = n \Leftrightarrow x\in D_n.$$
  It is easy to see that  $X$ is proper, but on the other hand $X$ does not have bounded
  geometry since the number of elements in $D_d(x_n, 1)$ tends to infinity
  for $n\rightarrow\infty$.
\end{ex}

\section{Construction of a \plig on $G$. }

In this and the section,  we will prove that every \lcsc group has a
 \plig  $d$. Together with theorem \bref{thm:plig_coarse_eq} we have
 therefore that a \lcsc group has exactly one coarse equivalence class of  \pligs.


%



\begin{df}
  We will call a  map $l:G\rightarrow \bR_+$  for a   {\em
  length function}\index{length function}, if it satisfies the following conditions \bref{eq:wbh1}--\bref{eq:wbh3}.

  \begin{eqnarray}\label{eq:wbh1}
    \forall_{g\in G}\quad
    l(g) = 0 &\Leftrightarrow& g = e.\\\label{eq:wbh2}
   \forall_{g\in G}\quad l(g)& = &l(g^{-1}).\\
 \label{eq:wbh3}
 \forall_{g, h\in G}\quad l(g\cdot h) &\leq& l(g)+l(h).
 \end{eqnarray}

 \end{df}

 \begin{lem}\label{lem:length_cond}
   \begin{enumerate}
     \item If $l:G\rightarrow \bR_+$ is a length function, then 
       \begin{equation}
         \label{eq:u11}
        d(x,y) = l(y^{-1}x), \quad x,y\in G   
       \end{equation}
       is a left invariant metric on $G$. 
\item
Conversily, if $d: G\times G\rightarrow \bR_+$ is a left invariant metric
on $G$, then 
$$l(x) = d(x,e), \quad x\in G$$
is a length function on $G$, and $d$ is the metric obtained from $l$ by \bref{eq:u11}.
   \end{enumerate}

Moreover, if $l$ is a length function on $G$ and $d(x,y) = l(y^{-1}x)$ is
the associated left invariant metric, then $d$ generates the topology on
$G$ if and only if 
\begin{equation}
   \label{eq:wbh5}
    \{l^{-1}[0, r):\quad r>0\}\mbox{  is a basis for the neighborhoods of }\ e\in G_0.
\end{equation}
Moreover, $d$ is proper if and only if:
\begin{equation}
   \label{eq:wbh4}
   \forall_{r>0}\quad l^{-1}([0,r])\quad \mbox{ is compact.}
\end{equation}
 \end{lem}

 \begin{proof}
   The proof is elementary, and will be left to the reader.
 \end{proof}

 \begin{rem}\label{rem:4.3.5}
   If $l$ is a length function on $G$ satisfying \bref{eq:wbh5}, then the
   associated metric 
   \begin{equation*}
     d(x,y) = l( y^{-1}x), \quad x,y \in G
   \end{equation*}
generates the given topology on $G$. Therefore 
\begin{equation*}
   l^{-1}([0,r]) = \{ x\in G:\quad d(x,e)\leq r  \}
\end{equation*}
is closed in $G$. Hence if we assume  \bref{eq:wbh5}, then  \bref{eq:wbh4}
is equivalent to that $l^{-1}([0,r])$ is relatively compact for all $r>0$,
which again is equialent to that 
\begin{equation*}
  B_d(e,n) = l^{-1}([0,n))
\end{equation*}
is relatively compact for all $n\in \bN$.
 \end{rem}

\begin{lem}\label{lem:ball_growth_limit}
   Let $G$ be a \lcsc group. Assume that the topology on $G$ is generated
   by a left invariant metric $\delta$, for which
   \begin{eqnarray}
     U = B_\delta(e,1)&\mbox{is relatively compact}\label{eq:cptl_gen}\\
     G = \bigcup^\infty_{k=1} U^k\label{eq:ass_cpct_gen}
   \end{eqnarray}

Then $G$ admits a left invariant metric $d$ generating the topology on $G$,
for which $B_d(e,1) =U$, and 
       \begin{equation}
         \label{eq:ball_growth_limit}
         \forall_{n\in\bN}\quad B(e,n) \subset B(e,1)^{2n-1}.
       \end{equation}

Moreover $d$ is a \plig on $G$.
\end{lem}

\begin{proof}
  Let $l_\delta(g) = \delta(g,e)$ be the length function associated to
  $\delta$, and define a function $l:G\rightarrow\bR_+$ by 
  \begin{equation}
     l(g) = \inf\bigg\{ 
\sum_{i=1}^k l_\delta(g_i):\quad  g = g_1, \ldots, g_k, \mbox{
  where }g_i\in U, \quad i =1, \ldots k, k\in\bN
\bigg\}
  \end{equation}

Clearly, $l(g)\geq 0$ for all $g\in G$, and by the assumption \bref{eq:ass_cpct_gen}, we have
that $l(g)<\infty$. Moreover, one checks easily that 
\begin{equation}
  l(g\cdot h) \leq l(g) + l(h), \quad g,h\in G.
\end{equation}

Since $U = U^{-1}$, we have also that
\begin{equation}
  l(g^{-1}) = l(g), \quad g\in G.
\end{equation}

Moreover
\begin{equation}\label{eq:uffe19}
  l(g) \geq l_\delta(g),\quad g\in G
\end{equation}
and
\begin{equation}\label{eq:uffe20}
  l(g) = l_\delta(g),\quad g\in U
\end{equation}

By \bref{eq:uffe19} and \bref{eq:uffe20}, we have that 
$$ l(g) = 0 \Leftrightarrow g= e.$$

Hence by lemma \bref{lem:length_cond}, $l$ is a length function on $G$. Let 
\begin{equation}\label{eq:uffe21}
  d(g,h) = l(g^{-1}\cdot h),\quad g,h\in G
\end{equation}
be the associated left invariant metric on $G$. By \bref{eq:uffe20} and
\bref{eq:uffe21} we have that
\begin{equation}
  \label{eq:uffe22}
  B_d(x,r) = B_\delta(x,r),\quad r\leq 1
\end{equation}
for all $g\in G$. Hence $d$ generates the same topology on $G$ as $\delta$
does. Moreover, we have that
\begin{equation}
  B_d(e,1) = B_\delta(e,1) =U.
\end{equation}

We next turn to the proof of \bref{eq:ball_growth_limit}. Let $n\in\bN$, and
let $g\in B_d(e,n)$. Then there exists a $k\in\bN$, and $g_1, \ldots,
g_k\in U$ such that 
$$ g = g_1\ldots g_k\quad\mbox{ and }\quad\sum_{i=1}^kl_\delta(g_i) < n.$$
Further we may assume, that $k\in\bN$ is minimal among all the numbers for
which such a representation is possible.

We claim, that in this case
\begin{equation}
  \label{eq:uffe23}
  l_\delta(g_i) + l_\delta(g_{i+1}) \geq 1,\quad i=1,\ldots, k-1.
\end{equation}

Assume namely, that $  l_\delta(g_i) + l_\delta(g_{i+1}) < 1$ for some $i$,
where $1\leq i\leq k-1$, then we have 
\begin{equation}
  l_\delta (g_i\cdot g_{i+1}) \leq   l_\delta(g_i) + l_\delta(g_{i+1}) < 1,
\end{equation}
and thus we have that $g_i\cdot g_{i+1}\in U$. Hence $g$ can be written as
a product of $k-1$ elements from $U$:

\begin{equation}
  g =  g_1\ldots g_i\cdot g_{i-1}(g_i\cdot g_{i+1})g_{i+2}\ldots g_k 
\end{equation}
for which 
\begin{equation}
  \sum_{j=1}^{i-1} l_\delta(g_j) + l_\delta( (g_i\cdot g_{i+1}) +
\sum_{j=i+2}^{k} l_\delta(g_j)  \leq 
  \sum_{j=1}^k l_\delta(g_j) < n.
\end{equation}

This contradicts the minimality of $k$, and hence \bref{eq:uffe23} must
hold. 

Let $\lfloor r\rfloor$ as usual denote the largest integer such that
\begin{equation*}
  \lfloor r\rfloor \leq r
\end{equation*}
>From \bref{eq:uffe23} we get that 

\begin{equation}
  \lfloor \frac{k}{2}\rfloor \leq 
\sum^{\lfloor \frac{k}{2}\rfloor}_{j=1}( l_\delta(g_{2j-1}) +
l_\delta(g_{2j}))
\leq
 \sum_{j=1}^k l_\delta(g_j) < n 
\end{equation}

Since both $\lfloor \frac{k}{2}\rfloor$ and $n$ are integers, we have that
$ \lfloor \frac{k}{2}\rfloor \leq n-1$, and therefore 
\begin{equation}
  k \leq 2\cdot \lfloor \frac{k}{2}\rfloor+1 \leq 2n-1,
\end{equation}
and hence we get that 
\begin{equation}
  g\in U^k \subset U^{2n-1} = B_d(e,1)^{2n-1},
\end{equation}
which prooves \bref{eq:ball_growth_limit}.

Since $ U^{2n-1} \subset  \overline{U}^{2n-1}$, where the latter set is
compact by assumption \bref{eq:cptl_gen}, we have that $B_d(e,n)$ is relatively
compact for all $n\in\bN$. Hence, by lemma \bref{lem:length_cond} and remark \bref{rem:4.3.5} $d$ is a \plig on $G$.


\end{proof}

We are now ready toprove, that there exists a \plig on every \lcsc group.



\begin{thm}\label{thm:plig_on_lcsc}
 Every \lcsc group $G$ has a \plig $d$. 
\end{thm}

\begin{proof}
   Let $G$ be \lcsc group. By Remark \cite[Theorem 1.22,page 34]{willard}, we can choose a left
   invariant metric $\delta_0$ on $G$, which generates the
   topology on $G$. Moreover, since $G$ is locally compact,there exists an
   $r>0$ such that the open ball $B_{\delta_0}(e,r)$ is relatively
   compact. Put now $\delta = \frac{1}{r}\delta_0$. Then $\delta$ is again
   a left invariant metric on $G$, which generates the topology. Moreower, 
   \begin{equation}
     U = B_\delta(e,1)
   \end{equation}
is relatively compact.

Put 
\begin{equation}
  G_0 = \bigcup^\infty_{k=1} U^k.
\end{equation}

Then $G_0$ is an open and closed subgroup of $G$. Since $G$ is second
countable, it follows that the space $Y = G/G_0$ of left $G_0$-cosets in
$G$ is a countable set. In the following we will assume that $|Y| =
\infty$. The proof in the case $|Y| < \infty$ can be obtained by the same
method with elementary modifications. 

We can choose a sequence $\{ x_n\}^\infty_{n=1} \subset G$, such that
$x_0=e$ and such that is a disjoint union of cosets:
\begin{equation}
  G = \bigcup^\infty_{n=0} x_n\cdot G_0.
\end{equation}

By lemma \bref{lem:ball_growth_limit}, there is a \plig
$d_0$ on $G_0$, such that 
\begin{eqnarray}
  \label{eq:uffe25}
  B_{d_0}(e,1) &=& U\\
  B_{d_0}(e,n)& \subset & U^{2n-1}.\label{eq:uffe26} 
\end{eqnarray}

In particular, $d_0$ is proper. Let 
\begin{equation}
  \label{eq:uffe27}
  l_0(h) = d_0(h,e),\quad h\in G_0
\end{equation}
be the length function associated with $d_0$. 

Put
\begin{equation}
  S = \{  x_1, x_2, x_3 \ldots \} \subset G
\end{equation}
and define $l_1:S\rightarrow \bN$ by
\begin{equation}
  \label{eq:uffe27,5}
  l_1(x_n) = n.
\end{equation}

Define furthermore a function $\tilde{l}: G\rightarrow [0, \infty [$ by setting
\begin{equation}
  \label{eq:uffe28}
  \tilde{l}(g) = \inf\bigg\{
l_0(h_0) + \sum^{k}_{i=1} (l_1(s_i) + l_0(h_i))
 \bigg\},
\end{equation}
where the infimum is taken over all the representations of $g$ of the form 
\begin{equation}
  \label{eq:uffe29}
  \begin{cases}
    g = h_0\cdot s_1 \cdot h_1 \cdot s_2\cdot h_2\cdots h_k\\
    k\in \bN\cup\{ 0  \}, h_0, \ldots, h_k\in G_0\\
    s_1, \ldots, s_k\in S
  \end{cases}
\end{equation}

Note that 
\begin{equation}
  G = G_0 \cup \bigcup^\infty_{n=1} x_n\cdot G_0 = G_0\cup S\cdot G_0
  \subset G_0 \cup G_0\cdot S\cdot G_0,
\end{equation}
so that every $g\in G$ has a representation of the form \bref{eq:uffe29}
with $k=0$ or $k=1$. 

Next, put
\begin{equation}
  \label{eq:uffe30}
  l(g) = \max \bigg\{ 
\tilde{l}(g), \tilde{l}(g^{-1})
\bigg\}, \quad g\in G
\end{equation}

We will show, that $l(g)$ is a length function on $G$, and that the
associated metric 
\begin{equation}
  \label{eq:uffe31}
  d(g,h) = l(g^{-1}\cdot h),\quad g,h\in G
\end{equation}
is a \plig on $G$.

It is easily checked that 
\begin{equation}
  \tilde{l}(g\cdot h) \leq \tilde{l}(g)+ \tilde{l}(h),\quad g,h\in G
\end{equation}

and hence also
\begin{equation}
  \label{eq:uffe32}
  l(g\cdot h) \leq l(g) + l(h), \quad g,h\in G.
\end{equation}

Moreover, by \bref{eq:uffe29}
\begin{equation}
  \label{eq:uffe33}
  l(g^{-1}) = l(g), \quad g\in G.
\end{equation}

If $g\in G$ and $l(g) < 1$, then $\tilde{l}(g) \leq l(h) < 1$. Thus for
    every $\epsilon>0$, $g$ has a representation of the form
    \bref{eq:uffe29}, such that
    \begin{equation}
      l_0(h_0) + \sum^k_{i=1} (l_1(s_i) + l_0(h_i)) < l(g) + \epsilon
    \end{equation}
and for sufficiently small $\epsilon$, we have that 
$$ l(g)+\epsilon < 1,$$
which implies that $k=0$, because 
$$ \forall_{s\in S}\quad  l_1(s) \geq 1.$$

Hence $g=h_0 \in G_0$, and 
\begin{equation*}
  l_0(g) = l_0(h_0) < l(g) +\epsilon.
\end{equation*}

Since $\epsilon$ can be chosen arbitrarily small, we have shown that
\begin{equation}
  \label{eq:uffe34}
  g\in G\mbox{ and }l(g)<1 \Rightarrow g\in G_0 \mbox{ and } l_0(g) \leq l(g).
\end{equation}

In particular, $l(g) =0 $ implies that $g=e$, which together with
\bref{eq:uffe32} and \bref{eq:uffe33} shows, that $l$ is a length function on
$G$. Hence by lemma \bref{lem:length_cond}
\begin{equation*}
  d(g,h) = l(g^{-1}\cdot h), \quad g,h\in G
\end{equation*}
is a left invariant metric on $G$. 

>From \bref{eq:uffe34} we have
\begin{equation}
  \label{eq:uffe35}
  g\in B_d(e,1) \Rightarrow g\in G_0\quad\mbox{ and }\quad l_0(g)\leq l(g).
\end{equation}

Conversely, if $g\in G_0$, then using \bref{eq:uffe28} with $k=0$ and
$h_0=g$, we get that 
$$ \tilde{l}(g) \leq l_0(g) $$
and therefore 
\begin{equation}
  \label{eq:uffe36}
  l(g) \leq \max \bigg\{
l_0(g), l_0(g^{-1})
\bigg\} = l_0(g), \quad g\in G_0.
\end{equation}

By \bref{eq:uffe35} and \bref{eq:uffe36}, we have 
\begin{equation}\label{eq:uffe36,5}
  B_d(e,r) = B_{d_0}(e,r), \quad 0<r\leq 1,
\end{equation}
and since $G_0$ is open in $G$, the sets
\begin{equation*}
  B_{d_0}(e,r), \quad 0<r\leq 1
\end{equation*}
form a basis of neighbourhoods for $e$ in $G$. Hence $d$ generates the
original topology on $G$.

It remains to be proved, that $d$ is proper, ie. that $B_d(e,r)$ is
relatively compact for all $r>0$. Note that it is sufficient to consider
the case, where $r= n \in \bN$.

Let $n\in \bN$, and let $g\in B_d(e,n)$. Then we have that 
$$ \tilde{l}(g) \leq l(g) < n. $$

Hence by \bref{eq:uffe29}, we see that
\begin{equation}
  \label{eq:uffe37}
 g = h_0\cdot s_1 \cdot h_1 \cdot s_2\cdot h_2\cdots s_k\cdot h_k,  
\end{equation}
where 
\begin{equation}
  \begin{cases}
    k\in \bN\cup\{ 0 \}\\
    h_0, \ldots, h_k\in G_0\\
    s_1, \ldots, s_k\in S\\
     l_0(h_0) + \sum^k_{i=1} (l_1(s_i) + l_0(h_i)) < n \label{eq:uffe38}
  \end{cases}
\end{equation}

Since 
$$l_1(s) \geq 1\quad \forall_{s\in S},$$ 
we have that $k\leq n-1$. Moreover,
since $l_1:S\rightarrow \bN$ is defined by 
\begin{equation*}
  l_1(x_m) = m,\quad m=1, 2, \ldots
\end{equation*}
we have
\begin{equation}
  \label{eq:uffe39}
  s_i \in \{ x_1, x_2, \ldots, x_{n-1} \}, \quad 1\leq i \leq k.
\end{equation}

Moreover
\begin{equation}
  \label{eq:uffe40}
  h_i \in B_{d_0}(e,n), \quad 0 \leq i \leq k
\end{equation}
because $l_0(h_i) < n $ by \bref{eq:uffe38}. Put
\begin{equation*}
  T(n) = \{ x_1, \ldots, x_{n-1} \} \cup \{ e \}
\end{equation*}

Then by \bref{eq:uffe37}, \bref{eq:uffe38}, \bref{eq:uffe39} and
\bref{eq:uffe40} we have
\begin{equation*}
  g \in B_{d_0}(e,n)\bigg( T(n)\cdot B_{d_0}(e,n)\bigg)^k
\subset \bigg( T(n)\cdot B_{d_0}(e,n)\bigg)^{k+1} 
\subset \bigg( T(n)\cdot B_{d_0}(e,n)\bigg)^n
\end{equation*}

where the last inclusion follows from the inequality $k\leq n-1$. Since
$g\in B_d(e,n)$ was chosen arbitrarily, we have shown that 
\begin{equation*}
  B_d(e,n) \subset \bigg( T(n)\cdot B_{d_0}(e,n)\bigg)^n.
\end{equation*}

But $d_0$ is a proper metric on $G_0$, and since $T(n)$ is a finite set, it
follows, that 
\begin{equation*}
  \bigg( T(n)\cdot \overline{B_{d_0}(e,n)}\bigg)^n
\end{equation*}
is compact.

Hence $B_d(e,n)$ is relatively compact for all $n\in \bN$, and therefore
$d$ is a proper metric on $G$, cf. remark \bref{rem:4.3.5}.

\end{proof}

\begin{rem}
  As mentioned in the introduction, Theorem \bref{thm:plig_on_lcsc} has previously been
  obtained in two important special cases, the compactly generated case
  \cite{lubotzky_etal} and the countable, discrete case \cite{tu}.
\end{rem}

\begin{ex}\index{\glnr}\index{length function!on \glnr}
In this example, we give an explicit formula for a \plig on \glnr. The same
formula will also define a \plig on every closed subgroup of \glnr. Recall
that
\begin{equation}
  \glnr = \{ A \in M_n(\bR): \det(A)\neq 0 \} 
\end{equation}
and the topology on \glnr is inherited from the topology of $M_n(\bR)\simeq
\bR^{n^2}$. We equip $M_n(\bR)$ with the operator norm 
\begin{equation*}
  ||A|| = \sup \{ ||Ax|| : x\in \bR^n, ||x||  \leq 1  \}, 
\end{equation*}
where $||x|| = \sqrt{x_1^2 + \cdots + x^2_n }$ is the Euclidian norm on
$\bR^n$. 

 Define a function on $GL(n, \bR)$ by
  \begin{equation}
    \label{eq:lgth_on_gl}
    l(A) =\max \{  \ln(1+ ||A-I||), \ln(1+||A^{-1}-I|| \}.
  \end{equation}

We claim that $l$ is a length function on \glnr, and that the associated
metric
$$ d(A, B) = l(B^{-1}A), \quad A, B\in \glnr$$
is a \plig on \glnr. We prove first, that $l$
is a length function. Clearly, $l(A) = l(A^{-1})$ and $l(A)= 0
\Leftrightarrow A = I$. 

Let $A, B \in \glnr$. Then  
$$ || A-I|| \leq e^{l(A)} -1, ||B-I||\leq e^{l(B)}-1.$$

Put $ X =  A- I$ and $Y = B-I$. Then
\begin{multline*}
|| AB -I || = || XY + X + Y || \leq ||X||\cdot||Y||+||X||+||Y||\\ =
   (||X||+1)(||Y||+1)-1 
   \leq e^{l(A)} e^{l(B)} -1
\end{multline*}
and hence
\begin{equation}\label{eq:u60}
  \ln(1+ (|| AB -I ||) \leq l(A)+ l(B).
\end{equation}

Substituting $(A, B)$ with $(B^{-1}, A^{-1})$ in this inequality, we get 
  \begin{equation}\label{eq:u61}
  \ln(1+ (|| (AB)^{-1} -I ||) \leq l(B^{-1}) + l(A^{-1}) = l(A) + l(B), 
\end{equation}
and by \bref{eq:u60} and \bref{eq:u61} it follows that $l(A+B) \leq
l(A)+l(B)$. Hence $l$ is a length function on \glnr.

To prove that $d$ is a \plig on \glnr, it suffices to check, that the
conditions \bref{eq:wbh4} and \bref{eq:wbh5} in lemma \bref{lem:length_cond} are fullfilled.

Since $A \rightarrow A^{-1}$ is a homeomorphism of \glnr onto itself,
\bref{eq:wbh5} is clearly fullfilled. To prove \bref{eq:wbh4} we let $r\in (0,
\infty)$ and put $M = e^r$. Since $||C|| \leq 1+ || C-I||$ for $C\in\glnr$,
we see that $l^{-1}([0,r])$ is a closed subset of 
$$ K = \{ A\in\glnr: ||A||\leq M, ||A^{-1}||\leq M  \}. $$

Denote
\begin{equation*}
  L = \{ (A,B)\in M_n(\bR)\times M_n(\bR): AB = BA = I, ||A||\leq M, ||B||\leq M  \},
\end{equation*}
then $L$ is a compact subset of $M_n(\bR)^2$, and $K$ is the range of $L$
by the continuous map $\pi: (A, B) \rightarrow A$ of $M_n(\bR)^2$ onto
$M_n(\bR)$. Hence $K$ is compact, and therefore $l^{-1}([0,r])$ is also
compact. This proves \bref{eq:wbh4}, and therefore $d$ is a \plig on \glnr.

\end{ex}

\begin{ex}
  Let $G$ be a connected Lie group. Then we can choose a left invariant
  Riemannian structure on $G$. Let $(g_p)_{p\in G}$ denote the
  corresponding inner product on the spaces $(T_p)_{p\in G}$. The path length metric on $G$
  corresponding to the Riemannian structure is
$$ d(g,h) = \inf_\gamma L(\gamma) $$
where
\begin{equation*}
  L(\gamma) = \int^b_a \sqrt{g_{\gamma(t)}( \gamma'(t), \gamma'(t))}dt
\end{equation*}
is the path length of a piecewise smooth path $\gamma$ in $G$, and where the
infimum is taken over all such paths, that starts in $\gamma(a)=g$ and ends
in $\gamma(b)= h$. Then $d$ is a left invariant metric on $G$ which induces
the given topology on $G$, cf. \cite[p.51-52]{helgason:diff_geom}.

We claim that $d$ is a proper metric on $G$. To prove this, it is
sufficient to prove, that $B_d(e,r)$ is relatively compact for all
$r>0$. Let $r>0$, and let $g\in B(e,r)$. Then $e$ and $g$ can be connected
with a piecewise smooth path $\gamma$ of length $L(\gamma) <r$. 

Now $\gamma$ can be divided in two paths each of length
$\frac{1}{2}L(\gamma)$. Let $h$ denote the endpoint of the first path. Then 
\begin{equation*}
  d(e,h) \leq \frac{L(\gamma)}{2}\quad\mbox{and}\quad d(h,g)\leq \frac{L(\gamma)}{2}.
\end{equation*}

Hence $g = h(h^{-1}g)$, where $d(e,r) < \frac{r}{2}$ and $d(h^{-1}g, e) =
d(g,h) < \frac{r}{2}$. This shows, that 
\begin{equation*}
  B(e,r)\subset B(e,\frac{r}{2})^2
\end{equation*}
and hence
\begin{equation*}
  B(e,r) \subset B(e,r\cdot 2^{-k})^{2^k}
\end{equation*}
for all $k\in\bN$. Since $G$ is locally compact, we can choose a $k\in\bN$
such that $B(e, r2^{-k})$ is relatively compact. Hence $B(e,r)$ is
contained in the compact set 
\begin{equation*}
  \overline{ B(e, r2^{-k})}^{2^k}.
\end{equation*}

This shows that $d$ is proper, and therefore $d$ is a \plig on $G$.

\end{ex}

\section{Exponentially controlled growth of the $d$-balls}

\begin{df}
   Let $(G,d)$ be a \lcsc group with a \plig, and let $\mu$ denote the Haar
    measure on $G$. Then we say that  {\em
      \expgrow}\index{growth!exponentially controlled} if 
there exists constants $\alpha, \beta>0$, such that 
        \begin{equation}
          \label{eq:sphere_growth}
           \mu (B_d(e,n))\leq \beta e^{\alpha n}, \quad\forall_{n\in \bN}.
        \end{equation} 

Note, that if \bref{eq:sphere_growth} holds, then 
 \begin{equation}
     \label{eq:expgrow_general}
     \mu( B_d(e,r)) \leq \beta'\cdot e^{c_2\cdot r},\quad r\in [1, \infty),
 \end{equation}
 where $\beta' = \beta\cdot e^{\alpha}$. 
  \end{df}

We now turn to the problem of constructing a \plig on $G$, for which
\expgrow . We first prove the following simple combinatorial lemma:

\begin{lem}\label{lem:combinatorial}
  Let $n\in \bN$, and let $k\in \{ 1, \ldots, n \}$.

Put
\begin{equation}
  N_{n,k} = \bigg\{ 
(n_1, n_2, \ldots, n_k)\in \bN^k:\quad \sum^k_{i=1} n_i \leq n
\bigg\}
\end{equation}

Then the number of elements in $N_{n,k}$ is
\begin{equation*}
  |N_{n,k}| = \binom{n}{k}.
\end{equation*}
\end{lem}

\begin{proof}
 The map
  \begin{equation*}
    (n_1, \ldots, n_k) \rightarrow \{ n_1, n_1+n_2, \ldots, n_1+n_2+\cdots
    + n_k \}
  \end{equation*}
is a bijection from $N_{n,k}$ onto the set of subsets of $\{ 1, \ldots, n
\}$ with $k$ elements, and the latter set has of course $\binom{n}{k}$ elements.  
\end{proof}

Having established lemma \bref{lem:combinatorial}, we can now turn to giving
a proof of the main theorem of this section:

\begin{thm}\label{thm:lcsc_has_expgrow}
  Every \lcsc group $G$ has a \plig $d$, for which \expgrow.
\end{thm}

\begin{proof}

The result is obtained, by modifying the construction of a \plig on $G$
from the proof of theorem \bref{thm:plig_on_lcsc}.

Let $U, G_0, d_0, d$ and $S= \{ x_1, x_2, \ldots \}$ be as in the proof of
theorem \bref{thm:plig_on_lcsc}, and note that by \bref{eq:uffe36,5} we have
\begin{equation*}
  B_d(e,1) = B_{d_0}(e,1) = U.
\end{equation*}

For each $i\in \bN$ the set $\overline{U}\cdot x_i$ is compact in $G$, and
can therefore be covered by finitely many left translates of $U$:
\begin{equation}
  \label{eq:uffe41}
  U\cdot x_i \subset \overline{U}\cdot x_i \subset \bigcup^{p(i)}_{j=1}
  y_{i,j}\cdot U.
\end{equation}

Define $l^*_1:S\rightarrow [1,\infty[$ by
\begin{equation}
  \label{eq:uffe42}
  l_1^*(x_i) = i+ \log_2(p(i)), 
\end{equation}
and note that 
\begin{equation*}
l^*_1(x_i) \geq i = l_1(x_i), \quad x_i\in S,  
\end{equation*}
where $l_1:S\rightarrow \bN$ is the map defined in \bref{eq:uffe27,5}.

We will now repeat the construction of the left invariant metric $d$ in the
proof of theorem \bref{thm:plig_on_lcsc}, with $l_1$ replaced by $l^*_1$,
ie. we first define a function $\tilde{l}^*:G\rightarrow [0,\infty[$ by
\begin{equation}
  \label{eq:uffe43}
\tilde{l}^*(g) = \inf \bigg\{
l_0(h_0) + \sum^k_{i=1} (l^*_1(s_i) + l_0(h_i))
\bigg\}
\end{equation}
where the infimum is taken over all representations of $g$ of the form
\begin{equation}
  \label{eq:uffe44}
  \begin{cases}
    g = h_0\cdot s_1\cdot h_1 \cdot s_2\cdot h_2 \cdots s_k\cdot h_k\\
    k\in \bN\cup \{ 0 \}\\
    h_0, \ldots, h_k \in G_0\\
    s_1, \ldots, s_k \in S
  \end{cases}
\end{equation}

Next, we put
\begin{equation}
  \label{eq:uffe44.1}
  l^*(g) = \max \bigg\{ 
\tilde{l}^*(g), \tilde{l}^*(g^{-1})
\bigg\}, \quad g\in G
\end{equation}
and 
\begin{equation}
  \label{eq:uffe44.2}
  d^*(g,h) = l^*(g^{-1}\cdot h),\quad g,h\in G.
\end{equation}

Then, exactly as for the metric $d$ in the
proof of theorem \bref{thm:plig_on_lcsc} we get that $d^*$ is a left
invariant metric on $G$, which generates the given topology on $G$,  and
which satisfies the following 
\begin{equation*}
    B_{d^*}(e,1) = B_{d_0}(e,1) = U,
\end{equation*}
(cf. the proof of theorem \bref{thm:plig_on_lcsc}).

Moreover, since $l_1^* \geq l_1$, we have that $d^*\geq d$, so the
properness of $d$ implies, that $d^*$ is also proper.

Let $n\in \bN$. Since 
\begin{equation*}
l^*(g) \leq \tilde{l}^*(g), \quad g\in G,
\end{equation*}
we get from \bref{eq:uffe43} and \bref{eq:uffe44} that
\begin{equation*}
  B_{d^*}(e,n) \subset \bigg\{ 
g\in G:\quad \tilde{l}^*(g) < n
\bigg\}
\end{equation*}

Hence every $g\in  B_{d^*}(e,n) $ can be written on the form
\bref{eq:uffe44} with 
\begin{equation*}
  l_0(h_0) + \sum^k_{i=1} (l^*_1(s_i) + l_0(h_i)) < n.
\end{equation*}

Note, that since 
\begin{equation*}
l^*_1(s) \geq l_1(s) \geq 1, \quad s\in S,
\end{equation*}
we have that $k\leq n-1$. 

Choose next natural numbers $m_0, \ldots, m_k$, such that 
\begin{equation*}
  l_0(h_i) < m_i \leq l_0(h_i) +1, \quad i= 0, \ldots, k.
\end{equation*}

Then we have that
\begin{equation*}
  h_i \in  B_{d^*}(e,m_i), \quad i = 0, \ldots, k,
\end{equation*}
and
\begin{equation*}
      m_0 + \sum^k_{i=1} (l_1^*(s_i) + m_i) < n + (k+1) \leq 2n+1.
\end{equation*}

Hence
\begin{equation}
  \label{eq:uffe45}
   B_{d^*}(e,n) \subset 
\bigcup_M\quad B_{d_0}(e, m_o )\cdot x_{n_1}\cdot B_{d_0}(e, m_1 )\cdot x_{n_2}\cdots
x_{n_k}\cdot B_{d_0}(e, m_k),
\end{equation}
where $M$ is the set of tuples
\begin{equation}
  \label{eq:uffe46}
  \bigg(
k, m_0, m_1, \ldots, m_k, n_1, \ldots, n_k
\bigg)
\end{equation}
for which 
\begin{equation*}
  \begin{cases}
    k\in \{ 0, \ldots, n-1 \}\\
  n_1, \ldots, n_k, m_0, \ldots, m_k \in \bN\\
 \sum_{i=0}^k m_i + \sum_{i=1}^k l^*_1(x_{n_i}) < 2n +1
  \end{cases}
\end{equation*}

By \bref{eq:uffe42}, the latter condition can be rewritten as 
\begin{equation}
  \label{eq:uffe47}
 \sum_{i=0}^k m_i + \sum_{i=1}^k ( n_i + log_2(p(n_i)) < 2n +1,
\end{equation}
where $p(n_i)\in \bN$ are given by formula \bref{eq:uffe41}. Since $m_i, n_i
\in \bN$ and $p(n_i)\geq 1$, it follows that 
\begin{equation*}
   \sum_{i=0}^k m_i + \sum_{i=1}^k n_i \leq 2n.
\end{equation*}

Hence $M \subset \bigcup^{n-1}_{k=0} M_k$, where $M_k$ is the set of
$2k+1$-tuples $(m_0, \ldots, m_k, n_1, \ldots, n_k)$ of natural numbers for
which
\begin{equation*}
   \sum_{i=0}^k m_i + \sum_{i=1}^k n_i \leq 2n.
\end{equation*}

Therefore, by lemma \bref{lem:combinatorial} we have that $|M_k| =
\binom{2n}{2k+1}$, and thus 
\begin{equation}
  \label{eq:uffe48}
  |M| \leq \sum^{n-1}_{k=0}\binom{2k}{2k+1} \leq
  \sum^{2n}_{j=0}\binom{2n}{j} = 2^{2n}. 
\end{equation}

>From \bref{eq:uffe26}, we have 
\begin{equation*}
  B_{d_0}(e, m_i) \subset U^{2m_i-1} \subset U^{2m_i}.
\end{equation*}

Hence by \bref{eq:uffe45}, we have
\begin{equation}
  \label{eq:uffe49}
    B_{d^*}(e,n) \subset 
\bigcup_M\quad U^{2m_o}\cdot x_{n_1}\cdot U^{2m_1} \cdot x_{n_2}\cdots
x_{n_k}\cdot U^{2m_k},
\end{equation}
where $|M| \leq 2^{2n}$ and where \bref{eq:uffe47} holds for each $(k, m_0,
\ldots, m_k, n_1, \ldots, n_k) \in M$. Since $\overline{U}^2$ is compact,
it can be covered by finitely many left translates of $U$, ie. 
\begin{equation*}
  U^2 \subset \overline{U}^2 \subset \bigcup^q_{i=1} z_i\cdot U,\quad z_1,
  \ldots z_q \in G.
\end{equation*}

It now follows that for
every $k\in \bN$, the set $U^k$ can be covered by $q^{k-1}$ translates of
$U$, namely
\begin{equation}
  \label{eq:uffe50}
  U^k \subset \bigcup^q_{i_1 = \cdots = i_{k-1} = 1}\quad z_{i_1}\cdots
  z_{i_{k-1}}\cdot U, \quad z_1, \ldots, z_q\in G.
\end{equation}

We can now use \bref{eq:uffe41} and \bref{eq:uffe50} to control the Haar measure of
each of the sets 
\begin{equation}
  \label{eq:uffe51}
  U^{2m_o}\cdot x_{n_1}\cdot U^{2m_1} \cdot x_{n_2}\cdots
x_{n_k}\cdot U^{2m_k},
\end{equation}
from \bref{eq:uffe49}. By \bref{eq:uffe50} we see that $U^{2m_0}$ can be
covered by $q^{2m_0-1}$ left translations of $U$. 

Combined with
\bref{eq:uffe41}, we get that $U^{2m_0}\cdot x_{n_1}$ can be covered by
$q^{2m_0-1}\cdot p(n_1)$ left translates of $U$. Hence 
\begin{equation*}
U^{2m_0}x_{n_1}U^{2m_1} = \bigcup_{w\in A_1}w\cdot U^{2m_1+1},  
\end{equation*}
where $|A_1| \leq 2^{2m_0 -1}p(n_1)$. 

Again by \bref{eq:uffe50}, we see that   $U^{2m_1+1}$ can be
covered by $q^{2m_1}$ left translations of $U$, so altogether we see that
the set 
$$ U^{2m_0}x_{n_1}U^{2m_1}$$
can be covered by  $q^{2m_0 + 2m_1-1}p(n_1)$ left translations of
$U$. Continuing this procedure, we get that the set in \bref{eq:uffe51} can
be covered by 
\begin{equation*}
   q^{2m_0 + 2m_1 + \cdots + 2m_k -1}p(n_1)\cdot p(n_2)\cdots p(n_k) 
\end{equation*}
left translates of $U$, and hence the Haar measure of the set satisfies
that
\begin{multline*}
  \mu\bigg(
 U^{2m_o}\cdot x_{n_1}\cdot U^{2m_1} \cdot x_{n_2}\cdots
x_{n_k}\cdot U^{2m_k}
\bigg) \\
\leq
q^{2\cdot \sum^k_{i=0} m_i}\cdot \prod^k_{i=1} p(n_i)\mu(U)
\\
\leq
q^{2\cdot \sum^k_{i=0} m_i}\cdot 2^{ \sum^k_{i=1} \log_2(p(n_i))}\mu(U)
\end{multline*}

By \bref{eq:uffe47}, we have that 
\begin{equation*}
  \begin{cases}
     \sum^k_{i=0} m_i \leq 2n+1 \\
\sum^k_{i=1}   \log_2(p(n_i)) \leq 2n+1
  \end{cases}.
\end{equation*}

Hence, we have that 
\begin{equation*}
   \mu\bigg(
 U^{2m_o}\cdot x_{n_1}\cdot U^{2m_1} \cdot x_{n_2}\cdots
x_{n_k}\cdot U^{2m_k}
\bigg) \leq (2q^2)^{2n+1}\cdot \mu(U).
\end{equation*}

This holds for all tupples $(k, m_0, \ldots, m_k, n_1, \ldots, n_k) \in M$,
and since we have shown that $|M| \leq 2^{2n}$, we get from \bref{eq:uffe45}
that
\begin{equation}
  \mu( B_{d^*}(e,n)) \leq (4q^2)^{2n+1}\mu(U),  \quad n\in\bN,
\end{equation}
which shows that the $d^*$-balls have exponentially controlled growth.

\end{proof}

\begin{ex}\label{thm:discrete_countable_gp}

In this example, we will give a more direct proof of Theorem
\bref{thm:lcsc_has_expgrow} in the case of a countable discrete group
$\Gamma$.  If $\Gamma$ is finitely generated, it is elementary to check,
that the word length metric $d$ is proper and that the $d$-balls have
exponentially controlled growth, so we can assume that $\Gamma$ is
generated by an infite, symmetric set $S$, such that $e\notin S$.

We can write $S$ as a disjoint union 
\begin{equation*}
  S = \bigcup^\infty_{n=1} Z_n,
\end{equation*}
where each $Z_n$ is of the form $\{ x_n, x_n^{-1} \}$. Note that  $|Z_n|= 2$
if $x_n \neq x_n^{-1}$, and $|Z_n|= 1$  if  $x_n = x_n^{-1}$. Define
a function 
\begin{equation*}
  l_0: S\rightarrow \bN
\end{equation*}
by
\begin{equation*}
  l_0(x_n) =  l_0(x_n^{-1}) = n.
\end{equation*}

Next, define a function 
\begin{equation*}
  l: \Gamma\rightarrow \bN\cup\{ 0 \}
\end{equation*}
by
\begin{equation}  \label{eq:length_df}
  l(g) =
  \begin{cases}   
    \inf\{ \sum_{k=1}^n l(g_k)  \} & g\neq e  \\
    0& g= e
  \end{cases}
\end{equation}
where the infimum is taken over all representations of $g$ of the form 
\begin{equation*}
  g = g_1\cdots g_n, \quad g_i\in S, n\in \bN.
\end{equation*}
Then it is easy to check, that $l$ is a length function on $\Gamma$, and
since 
\begin{equation*}
  l(g) \geq 1\quad\mbox{for}\quad g\in\Gamma\setminus\{ e \}
\end{equation*}
the associated left invariant metric 
\begin{equation}
  d(g,h) = l(g^{-1}h), \quad g,h\in \Gamma
\end{equation}
generates the discrete topology on $\Gamma$. Put
\begin{equation}
  D(e,n) = \{ g\in \Gamma: d(g, e) \leq n  \}.
\end{equation}

We will next show, that 
\begin{equation}
  \label{eq:u76}
  |D(e,n) | \leq 3^n, \quad n\in\bN,
\end{equation}
which clearly implies, that \expgrow. In order to prove \bref{eq:u76}, we will
show by induction in $n\in\bN$, that the set 
\begin{equation}
   \partial(e, n) = \{ g\in \Gamma: d(g,e) = n  \}
\end{equation}
satisfies
\begin{equation}
  \label{eq:u77}
  |\partial(e,n)| \leq 2\cdot 3^{n-1}, \quad n\in\bN.
\end{equation}

Since $l_0(s)\geq 2$ for $s\in S\setminus\{ x_1, x_1^{-1} \}$, we have for
$g\in\Gamma$, that
\begin{equation}
  l(g) = 1 \Leftrightarrow g \in \{ x_1, x_1^{-1} \}.
\end{equation}

Hence
\begin{equation}
  | \partial(e, 1) | = |\{ x_1, x_1^{-1} \}| \leq 2
\end{equation}
which proves \bref{eq:u77} for $n=1$. Let now $n\geq 2$ and assume as
induction hyphotesis, that 
\begin{equation}
  \label{eq:u78}
    |\partial(e,i)| \leq 2\cdot 3^{i-1}, \quad i = 1, \ldots, n-1.
\end{equation}

We shall then show, that 
\begin{equation*}
    |\partial(e,n)| \leq  2\cdot 3^{n-1}.
\end{equation*}

We claim, that 
\begin{equation}
  \label{eq:u79}
  \partial(e,n) \subset \bigcup^n_{k=1} Z_k\cdot \partial(e,n-k)
\end{equation}
 To prove
\bref{eq:u78}, let $g\in \partial(e,n)$. Then there exists a $m\in \bN$
and $g_1, \ldots, g_m\in S$ such that 
\begin{equation*}
  g = g_1\cdots g_m
\end{equation*}
and
\begin{equation*}
  \sum^m_{i=1} l_0 (g_i) = n
\end{equation*}

 Put 
$k = l_0(g_1)$. Then $k\in\bN$ and $k\leq n$. 
Now
 $$g = g_1\cdot (g_2\cdots g_m ), $$
where 
\begin{eqnarray}
  \label{eq:u80}
  l(g_1) \leq l_0( g_1) = k\\
\label{eq:u81}
 l(g_2\cdots g_m) \leq \sum^m_{i=2} l_0(g_i) = n-k.
\end{eqnarray}

But, since 
\begin{equation*}
  n = l(g) \leq l(g_1) + l(g_2 \cdots g_m)
\end{equation*}
equality holds in both \bref{eq:u80} and \bref{eq:u81}. Hence $g_2\cdots
g_m\in \partial(e, n-k)$, which proves \bref{eq:u79}. By \bref{eq:u79} we have

\begin{equation*}
   |\partial(e, n)| \leq
\sum_{k=1}^{n} |Z_k|\cdot|\partial(e, n-k)| \leq  2\sum_{k=1}^{n}
|\partial(e, n-k)| =  2\sum_{i=1}^{n-1} |\partial(e, i)|,
\end{equation*}

Since  $|\partial(e,0)| = |\{
e\}| =1 $, we get by the induction hypothesis \bref{eq:u78}, that 
\begin{multline}
   |\partial(e, n)| \leq  \sum_{i=0}^{n-1} 2|\partial(e, i)|
 \leq  2(1+ \sum_{i=1}^{n-1} 2\cdot 3^{i-1})
=  2(1+3^{n-1}-1)) = 2\cdot 3^{n-1}.
\end{multline}
This completes the proof of the induction step. Hence \bref{eq:u77} holds
for all $n\in\bN$. Since $l$ only takes integer values, we have 
$$ D(e,n) =  \bigcup_{i=0}^n \partial(e, i).$$
Therefore
\begin{equation}
  |D(e, n)| =  \sum_{i=0}^n |\partial(e, i)| \leq  1 +  \sum_{i=1}^{n}
  2\cdot 3^{i-1}  = 3^n.
\end{equation}

This proves \bref{eq:u76}, and it follows that \expgrow.

\end{ex}

\section{Affine actions on Banach spaces.}

We have shown in theorem \ref{thm:plig_on_lcsc} that for any \lcsc group $G$ there
exists a \plig $d$, and  we have shown in theorem \ref{thm:lcsc_has_expgrow} that $d$ can be
chosen so that \expgrow. We will now construct an affine action of $G$ on the
reflexive seperable strictly convex Banach space 
$$\lp.$$

Gromov suggested in \cite{gromov_novikov}, that it is purposeful to
attack the Baum Connes Conjecture by considering proper affine isometric
actions on various Banach spaces. 

It was shown by N. Higson and G. Kasparov in \cite{higson_kasparov} that
the Baum-Connes Conjecture holds for discrete countable groups that admit a proper affine
isometric action on a Hilbert space. In particular, this holds for all discrete
amenable groups. Moreover, Yu proved in
\cite{yu_hyperbolic}, that a word hyperbolic group $\Gamma$ has a proper
affine action on the uniform convex Banach space $l^p(\Gamma\times \Gamma)$
for some $p\in [2, \infty)$. 

Therefore, it is interesting to study what kind of proper affine isometric
actions on Banach spaces a given \lcsc group admits.

\begin{df}
  The group of affine actions on $G$:\index{group!affine}
  Let $X$ be a normed vector space, then the affine group of $X$ is:
  $$ \Aff(X)  = \{ \phi:X\rightarrow X\quad : \phi(x) =  Ax+b; A\in GL(X), b\in X  \}.$$
  We say that {\em  $G$ has an affine action on $X$}\index{action!affine}, if there exists a
  group homomorphism of $G$ on $\Aff(X)$, ie:
  \begin{equation}
 \alpha:G\circlearrowright \Aff(X),\label{eq:affine_action}
\end{equation}
  such that
  \begin{equation}
  \forall_{g,h\in G}\quad\alpha(g\cdot h) = \alpha(g)\circ
 \alpha(h).\label{eq:cont_action}
\end{equation}
  Let $\pi_g$ denote the linear part of $\alpha(g)$, and denote the
  translation part by $b(g)$.  
We say that $\alpha(g)$ is isometric if
 the linear part $\pi_g$ is isometric, \index{action!isometric}
i.e:
  $$ \forall_{\xi\in X}\quad ||\pi_g\xi|| = ||\xi|| .$$
  
Moreover, we say that the action is proper, if\index{action!proper}
  $$ \forall_{\xi\in X}\quad \lim_{g\rightarrow\infty} ||\alpha(g)\xi|| = \infty.$$
\end{df}

\begin{rem}
  Since $\alpha$ is a homomorphism of the group $G$ into $\Aff(X)$, we have
  that:
  \begin{eqnarray}
    \label{eq:cocycle_cond}
    \forall_{\xi\in X}& \quad \alpha(st)\xi &=  \alpha(s)(\alpha(t)\xi)
    \Leftrightarrow\notag\\
    \forall_{\xi\in X}&\quad\pi_{st}\xi +b(st) &= \pi_s\pi_t\xi + \pi_sb(t)+b(s) \Leftrightarrow \notag\\
    &\pi_{st} = \pi_s\circ \pi_t&\mbox{ and }\quad b(st) =
    \pi_s(b(t)) + b(s).
  \end{eqnarray}
  The formula for $b(st)$ is called the {\em cocycle condition with respect
  to $\pi$}\index{cocycle condition}.
\end{rem}

And we also need to know what a strictly convex space is -- and we will use
the opportunity to define a uniformly convex space as well:

\begin{df}
  Let $X$ be a  normed vector space, denote the unit ball by $S_X$.
  \begin{enumerate}
  \item The following two conditions are equivalent (see
    \cite[Prop. 5.1.2]{megginson}).
    If $X$ satisfies any of them, it is called {\em strictly convex}\index{convex!strictly}.
    \begin{enumerate}
    \item
      \begin{equation}
        \label{eq:rotund_1}
        \forall_{x\neq y\in S_X, 1>t>0}\quad ||tx + (1-t)y|| <1.
      \end{equation}
    \item 
      \begin{equation}
        \label{eq:rotund}
        \forall_{x\neq y\in S_X}\quad ||\frac{1}{2}(x+y)|| < 1.
      \end{equation}
    \end{enumerate}
    \item $X$ is called {\em uniformly convex}\index{convex!uniformly}
 if
      \begin{equation}
        \label{eq:uniformly_rotund}
        \forall_{\epsilon>0}\exists_{\delta>0}\forall{ x,y\in
        S_X}\quad ||x-y|| \geq
        \epsilon\Rightarrow ||\frac{1}{2}(x+y) || \leq 1-\delta
      \end{equation}
  \end{enumerate}
\end{df}

\begin{rem}
  Every space that is uniformly convex is also strictly convex (see
  \cite[Proposition 5.2.6]{megginson}).
  Examples of uniformly convex Banach spaces include 
$$l_p, l_p^n, \quad \infty>p>1, n\geq 1 $$ 
(this follows from Milman-Pettis theorem, see \cite[Theorem
  5.2.15]{megginson}).

  A uniformly convex Banach space is necessarily reflexive (see \cite[Theorem ]{megginson}).
 There are spaces that are strictly convex, but not uniformly convex, and also
  spaces that are strictly convex and not reflexive.  An 
  example of a strictly convex but not uniformly convex Banach space is:
  $$\bigoplus_{i=1}^\infty l^n_{p_n},\quad\mbox{ where } p_n = 1+\frac{1}{n}$$
  (with $l^2$ norm on the direct sum).
\end{rem}

As an application of the construction of a \plig with \expgrow on a given
\lcsc group in  theorem \bref{thm:lcsc_has_expgrow}, we will construct a proper
isometric action on the Banach space 
$$\lp.$$

 In \cite{brown_guentner} a proper isometric action is
constructed for a discrete group $\Gamma$ into the Banach space
$\bigoplus_{p=1}^\infty L^{p_n}(G, d\mu),$ where $p_n$ is an unbounded
sequence. We have generalized this result as follows:

\begin{thm}\label{thm:affine_action}\index{action!proper affine isometric
    on a \lcsc group}
 Let $G$ be a \lcsc
 group, and let $\mu$ denote the Haar measure. Then there exists a proper affine isometric action $\alpha$ of $G$ on the
seperable, strictly convex, reflexive Banach space 
\begin{equation*}
  X = \lp.
\end{equation*}
\end{thm}

\begin{proof}

  Let $G$ be as in the statement of the theorem, then according to theorem 
\bref{thm:plig_on_lcsc}
  and 
  theorem~\bref{thm:lcsc_has_expgrow} we can choose a \plig $d$ on $G$ where
  \expgrow, ie.
  \begin{equation}\label{eq:uffe2.6enhalv}
    \exists_{\alpha>0}\quad \mu(B_d(e,n) \leq \beta\cdot e^{\alpha n}, 
  \end{equation}
for some constants $\alpha, \beta > 0$. We can without loss of generality
assume that $\beta\geq 1$.

Consider the functions $phi^n_g: G -> R$ given by:
\begin{equation}
  \label{eq:lip_fn}
  \phi^n_x(y) = 
  \begin{cases}
    1 - \frac{d(x,y)}{n}&\mbox{when }d(x,y)\leq n\\
    0 &\mbox{when } d(x,y)\geq n
  \end{cases}
\end{equation}

It is easy to check, that $phi^n_g$ is $\frac{1}{n}$-Lipschitz: 
\begin{equation}\label{eq:uffe2.8}
 |\phi^n_x(y)-\phi^n_x(z)|\leq \frac{d(y,z)}{n}.  
\end{equation}

Assume that $x\in B_d(e, \frac{n}{2})$:
\begin{equation}\label{eq:bound_by_charac}
  \phi^n_e(x) =   1 - \frac{d(x,e)}{n} \geq 1 - \frac{\frac{n}{2}}{n} =
  \frac{1}{2}\cdot 1_{ B_d(e, \frac{n}{2})}(x), 
\end{equation}

Let $\Cucb(G)$ denote the set of uniformly continuous bounded functions from
$G$ to $\bR$. Define $b^n:G\rightarrow \Cucb(G)$ by:
\begin{equation}
  \label{eq:b_n}
  b^n(g) = \lambda(g)\phi^n_e - \phi^n_e \Rightarrow b^n(g)(h) = \phi^n_e(g^{-1}h)-\phi^n_e(h).
\end{equation}


Since $d(g,e) = d(g^{-1},e)$, we have that 
\begin{equation*}
\phi^n_e(g) = \phi^n_e(g^{-1}),
\quad g\in G.  
\end{equation*}

Hence, by \bref{eq:b_n} and \bref{eq:uffe2.8} we have
\begin{equation}
  \label{uffe:2.11}
 |b^n(g)| \leq \big| \phi^n_e(g^{-1}h)-\phi^n_e(h) \big| =
 \big| \phi^n_e(h^{-1}g)-\phi^n_e(h^{-1}) \big|
 \leq
\frac{d(h^{-1}g,h^{-1}) }{n}  \leq \frac{d(e,g)}{n}
\end{equation}

Since $b_n(g) =0$, when $x\not\in B_d(e,n)\cup B_d(g,n)$, it follows that 
\begin{equation*}
 | b^n(g)| \leq \frac{d(e,g)}{n}\cdot  1_{ B_d(e, n)\cup B_d(g, n)} .
\end{equation*}

Hence $b^n(g) \in L^2(G, \mu)$ and
\begin{equation*}
   || b^n(g)||^{2n}_{2n} \leq \big( \frac{d(e,g)}{n}\big)^{2n} \big(\mu( B_d(e, n)) +
   \mu(B_d(g, n))\big).
\end{equation*}

Therefore, by \bref{eq:uffe2.6enhalv} and the left invariance of $\mu$, we
have that 
\begin{equation*}
   || b^n(g)||^{2n}_{2n} \leq \big( \frac{d(e,g)}{n}\big)^{2n}\cdot
   2\beta e^{\alpha n}.  
\end{equation*}

Using now, that $\beta\geq 1$ and $\sum^\infty_{n=1} \frac{1}{n^2} =
\frac{\pi^2}{6} < 2$, we get
\begin{equation*}
    \sum^\infty_{n=1} || b^n(g)||^{2}_{2n} \leq  \sum^\infty_{n=1}
    \frac{d(e,g)^2}{n^2}(2\beta)^{\frac{1}{n}}e^\alpha
\leq 4\beta e^\alpha d(g,e)^2
\end{equation*}

Hence 
$$b(g) = \oplus^\infty_{n=1} b^n(g) \in X$$
and
\begin{equation*}
  ||b(g)||_X \leq 2\sqrt{\beta} e^\frac{\alpha}{2} d(g,e).
\end{equation*}

Let $\tilde{\lambda}$ denote the left regular representation of $G$ on $X =
\lp$. Clearly $\tilde{\lambda}(g)$ is an isometry of $X$ for every $g\in G$.
We show next, that $b(g)$ fulfills the cocycle condition
\begin{equation}
  \label{eq:2.12}
  b(st) = \tilde{\lambda}(s)b(t) + b(s), \quad s,t\in G
\end{equation}
and \bref{eq:2.12} follows from 
\begin{multline*}
   b^n(st) = \lambda(st)\phi^n_e - \phi^n_e \\ 
=
 \lambda(s)(\lambda(t)\phi^n_e - \phi^n_e ) + (\lambda(s)\phi^n_e -
 \phi^n_e ) 
= \lambda(s)b^n(t) + b^n(s), \quad s,t\in G, 
\end{multline*}
for alle $n\in \bN$. By \bref{eq:2.12} we can define a continuous affine
action $\alpha$ of $G$ on $X$ by 
\begin{equation}
  \alpha(g)\xi = \tilde{\lambda}(g)\xi + b(g), \quad x\in X, g\in G.
\end{equation}

The last thing to show is that the action is metrically proper. For $\xi\in
X$ and $g\in G$, we have
\begin{equation*}
  || \alpha(g)\xi|| = || \tilde{\lambda}(g)\xi + b(g) || \geq 
 || b(g)|| - ||  \tilde{\lambda}(g)\xi|| = ||b(g)|| - ||\xi||.
\end{equation*}

Hence, we only have to check, that 
\begin{equation*}
  ||b(g)||\rightarrow \infty\quad\mbox{when}\quad d(g,e)\rightarrow\infty.
\end{equation*}

Let $g\in G$ and assume that $d(g,e)>2$. Moreover, let $N(g)\in \bN$ denote
the integer for which 
\begin{equation*}
  \frac{d(g,e)}{2} -1 \leq N(g) < \frac{d(g,e)}{2}.
\end{equation*}

For $ n = 1\ldots, N(g)$, we have that 
$$ d(g,e) > 2N(g)  \geq 2n.$$

Hence 
$$ \overline{B(e,n)} \cap \overline{B(g,n)} = \emptyset,$$
which implies that $\phi^n_e$ and $\phi^n_g$ have disjoint
supports. Therefore
\begin{equation*}
  || b^n(g)||^{2n}_{2n} = || \phi^n_g - \phi^n_e||^{2n}_{2n} = || \phi^n_g ||^{2n}_{2n} + ||
\phi^n_e  ||^{2n}_{2n} \geq || \phi^n_e ||^{2n}_{2n}
\end{equation*}

Since we have that
\begin{equation*}
  \phi^n_e \geq \frac{1}{2}\cdot 1_{B(e,\frac{n}{2})}
\end{equation*}
it follows that 
\begin{equation*}
   || b^n(g)||^{2n}_{2n} \geq 2^{-2n}\mu(B(e,\frac{n}{2})) \geq 2^{-2n}\mu(B(e,\frac{1}{2})).
\end{equation*}

Hence
\begin{multline}
  ||b(g)||^2 \geq \sum^{N(g)}_{n=1} || b^n(g) ||^{2}_{2n} \geq
\frac{1}{4} \sum^{N(g)}_{n=1} \mu(B(e, \frac{1}{2}))^{\frac{1}{n}} \\
 \geq
\frac{N(g)}{4}\cdot\min\big\{ \mu(B(e,\frac{1}{2})), 1 \big\}.\qquad
\end{multline}

Since 
\begin{equation*}
  N(g) \geq \frac{d(g,e)}{2}-1
\end{equation*}
it follows that 
\begin{equation*}
  || b(g) || \rightarrow \infty\mbox{ for }\quad d(g,e)\rightarrow\infty.
\end{equation*}

\end{proof}

\begin{cor}
   Let $G$ be a locally compact, 2-nd countable
 group. Then $G$ has a uniform embedding into the seperable, strictly
 convex Banach space
$$ \lp.$$
\end{cor}

\begin{proof}
We will show, that the map $b:G\rightarrow X$ constructed in the proof of
theorem \bref{thm:affine_action} is a uniform embedding. By the proof of theorem
\bref{thm:affine_action}, we have that 
\begin{equation}
  \label{eq:100}
  c_1\sqrt{d(g,e)} \leq ||b(g) ||_X \leq c_2 d(g,e),
\end{equation}
when $d(g,e) \geq c_3$ for some positive constants $c_1, c_2, c_3$. By the
cocycle condition \bref{eq:cocycle_cond}, we get for $g, h\in G$ that 
\begin{equation*}
  b(g) = b(h(h^{-1}g)) = \tilde{\lambda}(h)b(h^{-1}g) + b(h).
\end{equation*}

Hence
\begin{equation*}
  || b(g) - b(h) ||_X = || \tilde{\lambda}(h)b(h^{-1}g)||_X = ||b(h^{-1}g)||_X.
\end{equation*}

Since $d(h^{-1}g, e) = d(g,h)$, we obtain by applying \bref{eq:100} to
$b(h^{-1}g)$, that 
\begin{equation*}
   c_1\sqrt{d(g,h)} \leq ||b(g) - b(h)||_X \leq c_2 d(g,h),
\end{equation*}
when $d(g,h) \geq c_3$. This proves, that $b$ is a uniform embedding.

 \end{proof}

\bibliographystyle{alpha}
\bibliography{midtvej_references}


\end{document}